\documentclass[11pt]{article}
\usepackage{amsmath}
\font\tenmsb=msbm10
\font\sevenmsb=msbm7
\font\fivemsb=msbm5
\newfam\msbfam
\textfont\msbfam=\tenmsb
\scriptfont\msbfam=\sevenmsb
\scriptscriptfont\msbfam=\fivemsb

\def\mathbb#1{{\fam\msbfam\relax#1}}
\newcommand{\sel}{\setlength}
\sel{\textwidth}{390pt}
\sel{\oddsidemargin}{1cm}
\sel{\evensidemargin}{0cm}
\sel{\parskip}{0.20\baselineskip plus 1pt minus .5pt}
\sel{\headheight}{12.7pt}
\sel{\headsep}{2.5\baselineskip}
\sel{\topskip}{10.5pt}
\sel{\footskip}{29.5pt }
\sel{\listparindent}{9pt}
\sel{\topmargin}{6pt}
\sel{\pretolerance}{200}
\sel{\tolerance}{400}
\sel{\headheight}{0cm}
\sel{\brokenpenalty}{200}
\sel{\parindent}{17pt}
\sel{\textheight}{48\baselineskip}
\font\resu=cmb10 at 9pt
\newcommand{\al}{\alpha}
\newcommand{\io}{+\infty}
\newcommand{\lm}{\left}
\newcommand{\rt}{\right}
\newcommand{\te}{\theta}
\newcommand{\be}{\beta}

\newtheorem{theo}{\rm TH\'EOR\`EME}
\newtheorem{lem}{\rm LEMME}

\begin{document}
\noindent Th\'eorie des Nombres/{\sl Theory of Numbers}
\par
\begin{center}{\bf La fonction Z\^eta de Riemann prend une infinit\'e \\
de valeurs irrationnelles aux entiers impairs.}
\end{center}
\begin{center}{\sc Tanguy Rivoal}
\end{center}
\footnotesize
{\resu}
{{\bf R\'esum\'e -}
Nous montrons que la dimension de l'espace
vectoriel engendr\'e sur les rationnels par 1 et les n premi\`eres valeurs
de la fonction Z\^eta de Riemann aux entiers impairs cro\^{\i}t au moins
comme un multiple de log(n). Il en r\'esulte l'irrationalit\'e d'une
infinit\'e de valeurs de la fonction Z\^eta aux  entiers impairs.}
\begin{center}{\bf There are infinitely many irrational values of \\
the Riemann Zeta function at odd integers.}
\end{center}
{\resu} {{\bf Abstract -}
We provide a lower bound for the dimension of the vector
space spanned over the rationals by 1 and by the values of the Riemann Zeta
function at the first n odd integers. We prove that this dimension
increases at least like a constant times  log(n). As a consequence, the
Zeta function takes infinitely many irrational values at odd  integers.\\}

\vspace{0.1cm}
\small
\noindent{1  - INTRODUCTION. -}
\noindent Hormis l'irrationalit\'e de $\zeta(3)$, d\'emontr\'ee par
R.~Ap\'ery [1],  peu de r\'esultats sont connus sur la nature
arithm\'etique des nombres
$\zeta(2n+1)=\sum_{k\ge 1}1/k^{2n+1}$ ($n$ entier $\ge 1$).
\noindent Dans cette note, nous esquissons la d\'emonstration du
th\'eor\`eme suivant, dont d\'ecoule l'irrationalit\'e d'une infinit\'e de
$\zeta(2n+1)$ :
\begin{theo}
Pour tout $\varepsilon>0$, il existe un entier $N(\varepsilon)$ tel que si
$n>N(\varepsilon)$,
\end{theo}
$$\text{dim}_{\mathbb{Q}}\lm(\mathbb{Q}\,+\mathbb{Q}\;\zeta(3)
+\cdots+\mathbb{Q}\;\zeta(2n-1)+\mathbb{Q}\;\zeta(2n+1)\rt)\geq
\frac{(1-\varepsilon)}{1+\log(2)}\log(n).$$
\noindent La d\'emonstration s'inspire du travail de Nikishin [7] sur les
approximants de Pad\'e \mbox{de type I} des fonctions polylogarithmes
$L_n(z)=\sum_{k\ge 0}z^k/(k+1)^n$
(pour $z\in\mathbb{C}\;$, $|z|<1$). De fa\c{c}on plus
pr\'ecise, ayant  fix\'e des entiers $a$ et $b$ tels que
$1\le b\le a$, il d\'etermine, pour $|z|>1$,  des
polyn\^omes $Q_{i,n}(z)$ de degr\'e $\le n$ si
$i=1,\ldots,\, b$ et de degr\'e
$\le n-1$ si $i=0,\, b+1,\ldots,\, a$, tels que l'ordre en $z=\infty$ de la
fonction
$$
N_{n,a,b}(z)=Q_{0,n}(z)+\displaystyle\sum_{i=1}^a Q_{i,n}(z)L_i(1/z)
$$
soit au moins $an+b-1$ . En particulier, il obtient la formule explicite
$$
N_{n,a,b}(z)=\sum_{k=0}^{\io}
\frac{k(k-1)\cdots(k-an-b+2)}{(k+1)^a(k+2)^a\cdots(k+n)^a(k+n+1)^b}z^{-k}
$$
ce qui lui permet de montrer que
si $p/q\in\mathbb{Q}$ est tel que \mbox{$|q|>|p|^a(4a)^{a(a-1)}$,} alors les
nombres $1$, $L_1(p/q),\ldots,L_a(p/q)$ sont lin\'eairement ind\'ependants
sur $\mathbb{Q}$. Malheureusement les approximations de Nikishin,
sp\'ecialis\'ees en $z=-1$ et $b=a$, permettent seulement de montrer qu'il
y a au moins un irrationnel  parmi les nombres $\log(2)$, $\zeta(2)$,
$\zeta(3),\ldots,\zeta(a)$ (ce qui  r\'esulte a priori de la transcendance
de $\log(2)$, par exemple).\\

\noindent Pour am\'eliorer ce r\'esultat, on pourrait modifier la s\'erie
$N_{n,a,a}(z)$ en introduisant un param\`etre $r$ conduisant \`a de
meilleures estimations sur la croissance  des coefficients de la
combinaison lin\'eaire des valeurs de la fonction Z\^eta. Un choix
convenable de $r$ montrerait alors que la dimension $D(a)$ de l'espace
vectoriel  engendr\'e sur $\mathbb{Q}$ par $1,\zeta(2),
\zeta(3),\ldots,\zeta(a)$ est au moins
$c_0 \log(a)$ (o\`u $c_0$ est une constante effective). Cependant la formule
d'Euler $\zeta(2n)=2^{2n-1}B_n\pi^{2n}/(2n)!$ et la transcendance de $\pi$
impliquent que $D(a)\ge a/2$ : pour obtenir le Th\'eor\`eme 1, il s'agit
donc d'\'eliminer les nombres $\zeta(2n)$. Dans le cas de $\zeta(2),
\zeta(3)$ et $\zeta(4)$, K. Ball [1] a construit la s\'erie
$$
B_n=n!^2\sum_{k=1}^{\io}
\lm(k+\frac{n}{2}\rt)
\frac{(k-1)\cdots(k-n)(k+n+1)\cdots(k+2n)}{k^4(k+1)^4\cdots (k+n)^4}
$$
dont la forme particuli\`ere permet en effet d'\'eliminer $\zeta(2)$ et
$\zeta(4)$. Dans un message \`a l'auteur, K. Ball indiquait que sa
formule \'etait << facilement g\'en\'eralisable \`a $\zeta(5)$ et ainsi de
suite >> [2]. Nous g\'en\'eraliserons ici les s\'eries de Nikishin et Ball
en consid\'erant la  s\'erie (convergente pour $|z|\ge 1$)
\begin{eqnarray*}
S_n(z)&=&\sum_{k=0}^{\io} n!^{a-2r}
\frac{(k-rn+1)_{rn}(k+n+2)_{rn}}{(k+1)_{n+1}^a}z^{-k}
=\sum_{k=0}^{\io}R_{n}(k)z^{-k}
\end{eqnarray*}
o\`u $n$, $r$ et $a$ sont des entiers
v\'erifiant $1\leq r< a/2$, $n\in\mathbb{N}\,$
et o\`u $(\al)_k$ est le symbole de Pochammer : $(\al)_0=1$ et
$(\al)_k=\al(\al+1)\cdots(\al+k-1)$ si $k=1, 2,\ldots$.
\\ Ces s\'eries, sp\'ecialis\'ees en $z=1$, donneront des combinaisons
lin\'eaires \`a coefficients rationnels des Z\^eta impairs. Moyennant un
bon choix de $r$, ces combinaisons auront une  d\'ecroissance rapide vers
$0$ et leurs coefficients auront des d\'enominateurs et une croissance  bien
contr\^ol\'es. Le Th\'eor\`eme 1 d\'ecoulera alors du r\'esultat suivant,
d\^u \`a Y. Nesterenko [6] :\\

\noindent {\rm CRIT\`ERE D'IND\'EPENDANCE LIN\'EAIRE}
\textit{\noindent Consid\'erons $N$ r\'eels $\te_1,\te_2,\ldots,\te_N$
($N\geq 2$) et supposons qu'il existe $N$ suites d'entiers
$(p_{i,n})_{n\geq 0}$ tels que :
\begin{itemize}
\item[i)] $\log\left|\sum_{i=1}^N p_{i,n}\te_i\right|=
n\log(\al)+o(n)$ avec $0<\al<1$ ;
\item[ii)] $\forall i=1,\ldots,N$, $\log|p_{i,n}|\leq n\log(\be)+o(n)$ avec
$\be>1$.
\end{itemize}
Dans ces conditions, }
$\text{dim}_{\mathbb{Q}}(\mathbb{Q}\;\te_1+\mathbb{Q}\;
\te_2+\cdots+\mathbb{Q}\;\te_N)\geq 1-\log(\al)/\log(\be)$.\\

\noindent Je tiens \`a remercier vivement le professeur K. Ball : sans les
fructueux \'echanges que nous avons eus autour de sa s\'erie, cet article
n'aurait pu voir le jour. Je tiens \'egalement \`a exprimer toute ma
gratitude aux professeurs  F. Amoroso et M. Waldschmidt pour leurs
pr\'ecieux conseils et leur soutien constant.\\

\vspace{0.1cm}
\small
\noindent{2 - R\'ESULTATS AUXILIAIRES. -}
Pour $i=1,\ldots,a$, $j=0,\ldots,n$, d\'efinissons les nombres rationnels
$c_{i,j,n}=D_{a-i}\lm(R_{n}(t)(t+j+1)^a\rt)_{\vert t=-j-1}\;$
o\`u $D_{\lambda}=\frac{1}{\lambda!}d^{\lambda}/dt^{\lambda}$
et les polyn\^omes \\
$P_{0,n}(z)=-
\displaystyle\sum_{i=1}^a
\displaystyle\sum_{j=1}^n c_{i,j,n}
\displaystyle\sum_{k=0}^{j-1}\frac{1}{(k+1)^i} z^{j-k}$ et
$P_{i,n}(z)=\displaystyle\sum_{j=0}^n c_{i,j,n} z^j \quad(i=1,\ldots,a)$.\\
\begin{lem}
Si $n$ est pair et $a$ impair $\geq 3$, alors
$$S_n(1)=P_{0,n}(1)+\sum_{i=1}^{(a-1)/2} P_{2i+1,n}(1)\zeta(2i+1).
$$
\end{lem}
{\sl Preuve.} -
\noindent En d\'ecomposant $R_{n}(t)$ en fractions partielles, on a
$$R_{n}(t)= \displaystyle\sum_{i=1}^a\displaystyle\sum_{j=0}^n
\displaystyle\frac{c_{i,j,n}}{(t+j+1)^i}.$$
D'o\`u si $|z|>1$,
$$
S_{n}(z)=P_{0,n}(z)+\displaystyle\sum_{i=1}^a P_{i,n}(z)L_i(1/z).
$$
\noindent La convergence de la s\'erie $S_{n}(1)$ implique que
$\displaystyle \lim_{{z\to 1\atop |z|>1}}(P_{1,n}(z)L_1(1/z))=0$.
On peut \'ecrire
$c_{i,j,n}=(-1)^{a-i}D_{a-i}(\Phi_{n,j}(x))_{\vert x=j}$ o\`u
$\Phi_{n,j}(x)=R_{n}(-x-1)(j-x)^a$ : en appliquant l'identit\'e $(\al)_l
=(-1)^l(-\al-l+1)_l$, on montre que
$$\Phi_{n,n-j}(n-x)=(-1)^{na}\Phi_{n,j}(x).$$
Donc pour tout $k\geq 0$,
$\Phi_{n,n-j}^{(k)}(n-x)=(-1)^k(-1)^{na}\Phi_{n,j}^{(k)}(x)$~:
en particulier avec $k=a-i$ et $x=j$, on a $c_{i,n-j,n}=(-1)^{a-i}(-1)^{an}
c_{i,j,n}$, d'o\`u
$$
P_{i,n}(1)=(-1)^{(n+1)a+i}P_{i,n}(1).
$$
Si $n$ est pair et $a$ impair, on en d\'eduit que pour tout $i$ pair,
$P_{i,n}(1)=0$.\\

\noindent Le Lemme suivant donne une expression int\'egrale similaire \`a
celles de Beukers [4] (voir aussi [5], \S 1.3).
\begin{lem}
La s\'erie $S_{n}(z)$ admet la repr\'esentation int\'egrale,
pour $|z|\ge 1$~:
\begin{eqnarray*}
S_{n}(z)=\frac{((2r+1)n+1)!}{n!^{2r+1}z^{-(r+1)n-2}}
\int_{[0,1]^{a+1}}\lm(\frac{\prod_{i=1}^{a+1}x_i^{r}(1-x_i)}
{(z-x_1x_2\cdots x_{a+1})^{2r+1}}\rt)^n
\frac{dx_1dx_2\cdots dx_{a+1}}{(z-x_1x_2\cdots x_{a+1})^2}.
\end{eqnarray*}
\end{lem}
{\sl Preuve.} - Si $|z|>1$, cette \'egalit\'e s'obtient en d\'eveloppant en
s\'erie enti\`ere le \mbox{d\'enominateur} de la fraction sous le signe
int\'egral, l'interversion des signes  somme et int\'egral \'etant alors
justifi\'ee. En utilisant un argument de continuit\'e, on  montre que
l'\'egalit\'e reste \mbox{valable si $|z|=1$.}\\

\noindent Cette repr\'esentation int\'egrale permet alors d'estimer la
d\'ecroissance des nombres $S_n(1)$~:
\begin{lem}
La limite $s_{r,a}=\displaystyle \lim_{n\to\io}\lm
\vert S_{n}(1)\rt\vert^{1/n}$ existe et v\'erifie
\begin{eqnarray*}
s_{r,a}\le
(2r+1)^{2r+1}\frac{(ra+r)^{ra+r}(a-2r)^{a-2r}}{(ra+a-r)^{ra+a-r}}.
\end{eqnarray*}
\end{lem}
\noindent Pour estimer la croissance des nombres $P_{i,n}(1)$, il suffit de
majorer convenablement les coefficients $c_{l,j,n}$ au moyen de la formule
de Cauchy
$$
c_{l,j,n}=\frac{1}{2i\pi}\int_{\vert z+j+1\vert =1/2}R_n(z)(z+j+1)^{l-1}dz
$$
o\`u $\vert z+j+1\vert =1/2$ d\'esigne le cercle de centre $-j-1$ et 
de rayon $1/2$.
On obtient alors le
\begin{lem} Pour tout $i=0,\ldots,a$, on a
$\displaystyle\limsup_{n\to\io}\lm\vert P_{i,n}(1)\rt\vert^{1/n}
\leq 2^{a-2r}(2r+1)^{2r+1}.$
\end{lem}
Enfin, pour construire des combinaisons lin\'eaires \`a coefficients entiers,
il reste \`a d\'eterminer les d\'enominateurs des nombres $P_{i,n}(1)$, ce
qui r\'esulte du
\begin{lem} On pose $d_n=\text{ppcm}(1,2,\ldots,n)$. Alors pour
$i=0,\ldots,a$, $d_n^{a-i}P_{i,n}(1)\in\mathbb{Z}$.
\end{lem}
{\sl Preuve.} -
\noindent
Il s'agit d'\'evaluer le d\'enominateur commun des coefficients
$c_{i,j,n}$. Pour cela, fixons $n$ et $j$ et d\'ecomposons le num\'erateur
de $R_{n}(t)$ en $2r$ produits de $n$ facteurs cons\'ecutifs~:
on a $R_n(t)(t+j+1)^a=F_{1}(t)\cdots F_{r}(t) G_{1}(t)\cdots G_{r}(t)
H(t)^{a-2r}$
o\`u
$$
F_{l}(t)=\displaystyle\frac{(t-nl+1)_n}{(t+1)_{n+1}}(t+j+1)\,,\;
G_{l}(t)=\displaystyle\frac{(t+nl+2)_n}{(t+1)_{n+1}}(t+j+1)\,,\;
H(t)=\displaystyle\frac{n!(t+j+1)}{(t+1)_{n+1}}.
$$
En d\'ecomposant $F_{l}(t)$, $G_{l}(t)$  et $H(t)$ en fractions partielles,
on montre que pour tout entier $\lambda\geq 0$,
$d_n^{\lambda}(D_{\lambda}F_l)_{\vert t=-j-1}$,
$d_n^{\lambda}(D_{\lambda}G_l)_{\vert t=-j-1}$ et
$d_n^{\lambda}(D_{\lambda}H_l)_{\vert t=-j-1}$
sont des entiers. Gr\^ace \`a la formule de Leibniz, on en d\'eduit que
$d_n^{a-i}c_{i,j,n}\in\mathbb{Z}$ et donc 
$d_n^{a-i}P_{i,n}(1)\in\mathbb{Z}$ pour
$i=0,\ldots,a$ et pour tout $n\in\mathbb{N}$.\\

\vspace{0.1cm}
\small
\noindent{3 - PREUVE DU TH\'EOR\`EME 1. -}
\noindent Soit $a$ un entier impair $\geq 3$~: notons $\delta(a)$ la
dimension de l'espace vectoriel engendr\'e sur $\mathbb{Q}\;$ par $1$ et
les $\zeta(j)$ pour
$3\leq j\leq a$ et $j$ impair.\\
\noindent D'apr\`es le Th\'eor\`eme des Nombres Premiers,
$d_n=e^{n+o(n)}$. D\'efinissons pour tout entier $n\geq 0$~:
$\ell_{n}=d_{2n}^{a}S_{2n}(1)$, $p_{0,n}=d_{2n}^aP_{0,2n}(1)$ et
$ p_{i,n}=d_{2n}^aP_{2i+1,2n}(1)$ ($i=1,\ldots,(a-1)/2$).
Le Lemme~5 implique que $p_{i,n}\in\mathbb{Z}$ pour tout
$i=0,\ldots,(a-1)/2$ et d'apr\`es le \mbox{Lemme~1,}
$$
\ell_n=p_{0,n}+\displaystyle\sum_{i=1}^{(a-1)/2}p_{i,n}\zeta(2i+1).
$$
On peut appliquer le crit\`ere de Nesterenko avec $N=(a+1)/2$,
$\be=(e^a2^{a-2r}(2r+1)^{2r+1})^2$ (Lemme~4) et $\al=(e^as_{r,a})^2$
(Lemme~3)~:  pour tout entier $r$ tel que $1\leq r <a/2$, on en d\'eduit
alors $\delta(a)\geq f(a,r)/g(a,r)$ o\`u
\begin{eqnarray*}
f(a,r)&=& (a-2r)\log(2)+(ra+a-r)\log(ra+a-r)\\
&& -(ra+r)\log(ra+r)-(a-2r)\log(a-2r)
\end{eqnarray*}
et
$$
g(a,r)=a+(a-2r)\log(2)+(2r+1)\log(2r+1).
$$
Effectuons maintenant un d\'eveloppement limit\'e pour $a,r\to\io$ des
fonctions $f(a,r)$ et $g(a,r)$~:
$$
f(a,r)=a\log(r)+O(a)+O(r\log(r))\, \text{ et }\,
g(a,r)=(1+\log(2))a+O(r\log(r)).
$$
On choisit $r=r(a)$ comme l'entier $<a/2$ le plus proche de
$a(\log(a))^{-2}$~: on a alors $a\log(r)= a\log(a)(1+o(1))$ et
$r\log(r)=o(a)$. D'o\`u
$$
\delta(a)\geq \frac{f(a,r)}{g(a,r)}=
\frac{a\log(a)(1+o(1))+O(a)}{(1+\log(2))a+o(a)}=\frac{\log(a)}{1+\log(2)}
(1+o(1))\,,
$$
ce qui prouve le Th\'eor\`eme 1.\\

\vspace{0.2cm}
\centerline{\textbf{R\'ef\'erences}}
\vspace{0.2cm}

\noindent
$[1]$ R. Ap\'ery, \textit{Irrationalit\'e de $\zeta(2)$ et $\zeta(3)$},
Ast\'erisque \textbf{61}, 11-13 (1979).\\
$[2]$ K. Ball, Communication personnelle du 17 d\'ecembre 1999.\\
$[3]$ K. Ball, Communication personnelle du 4 janvier 2000.\\
$[4]$ F. Beukers, \textit{A note on the irrationality of $\zeta(2)$
and $\zeta(3)$},  Bull. London. Math. Soc. \textbf{11}, no. 33, 268-272
(1978).\\
$[5]$ R. Dvornicich et C. Viola, \textit{Some remarks on Beukers' integrals},
Number theory, Vol. II (Budapest, 1987), 637--657, Colloq. Math. Soc.
J\'anos Bolyai, 51, North-Holland, Amsterdam, 1990.\\
$[6]$ Y.V. Nesterenko, \textit{On the linear independence of numbers},
Mosc. Univ. Math. Bull. \textbf{40}, no. 1, 69-74 (1985) traduction de
Vest. Mosk. Univ., Ser. I, no. 1,  46-54 (1985).\\
$[7]$ E.M. Nikishin, \textit{On the irrationality of the values of the
functions $F(x,s)$}, Mat. Sbornik \textbf{37}, no. 3, 381-388 (1979).\\

\end{document}